\documentclass[12pt, nonatbib]{elsarticle}

\usepackage{amsmath}
\usepackage{mathscinet}
\usepackage{amssymb}
\usepackage{amsthm}
\usepackage{graphicx}
\usepackage[colorlinks]{hyperref}
\AtBeginDocument{%
  \hypersetup{%
    linkcolor=blue,%
    citecolor=green,%
  }%
}
\usepackage{cleveref}
\usepackage{geometry}
\makeatletter
\let\c@author\relax
\makeatother
\usepackage[backend=bibtex,style=numeric,sorting=nyt,
giveninits=true,isbn=false,url=false]{biblatex}
\renewbibmacro{in:}{\ifentrytype{article}{}{\printtext{\bibstring{in}\intitlepunct}}}
\DeclareFieldFormat*{title}{#1}
\crefname{equation}{}{}
\crefname{figure}{}{}

\newtheorem{theorem}{Theorem}[section]
\newtheorem{lemma}[theorem]{Lemma}

\theoremstyle{definition}
\newtheorem{definition}[theorem]{Definition}

\theoremstyle{remark}
\newtheorem{remark}[theorem]{Remark}

\numberwithin{equation}{section}

\journal{~~}

\begin{document}

\begin{frontmatter}

\title{Liouville theorems for fully nonlinear elliptic equations on half spaces}


\author{Yuanyuan Lian}
\ead{lianyuanyuan.hthk@gmail.com; yuanyuanlian@correo.ugr.es}

\address{Departamento de An\'{a}lisis Matem\'{a}tico,
Instituto de Matem\'{a}ticas IMAG, Universidad de Granada, Granada, Espa\~{n}a}

\tnotetext[t1]{This research is supported by the Grants PID2020-117868GB-I00 and PID2023-150727NB-I00 of the MICIN/AEI.}

\begin{abstract}
In this note, we prove two Liouville theorems for fully nonlinear uniformly elliptic equations on half spaces. The main tools are the boundary pointwise regularity, the Hopf type estimate and the Carleson type estimate. Our new proof is rather short.
\end{abstract}


\end{frontmatter}

\section{Introduction}
\label{intro}
In this note, we prove two Liouville theorems on half spaces for viscosity solutions of fully nonlinear uniformly elliptic equations
\begin{equation}\label{main}
F(D^2u)=0\quad\mbox{on}~~\mathbb{R}^n_+,
\end{equation}
where $F$ is uniformly elliptic (see \Cref{de.F}) and $\mathbb{R}^n_+$ is the standard half space ($n\geq 2$). Liouville-type theorems constitute one of the fundamental topics in the theory of partial differential equations. The first Liouville theorem was presented by Liouville in 1844 and immediately proved by Cauchy\cite{FARINA200761}. Large amount research has been performed for different types of elliptic equations (see \cite{FARINA200761} for a comprehensive survey on Liouville-type theorems).

For fully nonlinear uniformly elliptic equations, Braga \cite{MR3735564} gave a simple proof of a Liouville theorem (i.e., \Cref{t.ful-1}) on half-spaces. Armstrong, Sirakov and Smart \cite[Theorem 1.2]{MR2947535} extended it to general cones under an additional growth condition. In this note, we presents a new proof of this Liouville theorem on half spaces with the aid of the boundary pointwise regularity. In addition, with the same idea, we prove a second order Liouville theorem on half spaces (i.e., \Cref{t.ful-2}).

We use standard notations in this paper. For $x\in \mathbb{R}^n$, we may write $x=(x_1,...,x_n)=(x',x_n)$ and $|x|$ is the Euclidean norm of $x$ as usual. Let $ \mathbb{R}^n _+=\left\{x: x_n>0\right\}$ denote the upper half-space. Set $B_r(x_0)=\{x: |x-x_0|<r\}$ and $B_r=B_r(0)$. Additionally, the half ball is represented as $B_r^+(x_0)=B_r(x_0)\cap  \mathbb{R}^n _+$ and $B_r^+=B^+_r(0)$. Similarly, $T_r(x_0)\ =\{(x',x_{0,n}): |x'-x_0'|<r\}$ and $T_r=T_r(0)$ ($x_{0,n}$ is the $n$-th component of $x_0$). The $\{e_i\}_{i=1}^{n}$ is the standard basis of $ \mathbb{R}^n $.

In this note, the only assumption on $F$ is the uniform ellipticity. For the readers' convenience, we recall the definition below.
\begin{definition}\label{de.F}
The $F$ is called uniformly elliptic with ellipticity constants $0<\lambda\leq \Lambda$ if for any $M,N\in \mathcal{S}$,
\begin{equation}\label{e1.0}
\mathcal{M}^-(M,\lambda,\Lambda)\leq F(M+N)-F(N)
\leq \mathcal{M}^+(M,\lambda,\Lambda),
\end{equation}
where $\mathcal{S}$ denotes the set of $n\times n$ symmetric matrices. The $\mathcal{M}^-$ and $\mathcal{M}^+$ are Pucci's extremal operators, i.e.,
\begin{equation*}
\mathcal{M}^-(M,\lambda,\Lambda)=\lambda\sum_{\lambda_i>0} \lambda_i+
\Lambda\sum_{\lambda_i<0} \lambda_i, \quad
\mathcal{M}^+(M,\lambda,\Lambda)=\Lambda\sum_{\lambda_i>0} \lambda_i+
\lambda\sum_{\lambda_i<0} \lambda_i,
\end{equation*}
where $\lambda_i$ are eigenvalues of $M$.
\end{definition}

Note that Pucci's extremal operators $\mathcal{M}^-$ and $\mathcal{M}^+$ are uniformly elliptic. In this note, we consider viscosity solutions (see \cite{MR1351007} and \cite{MR1118699} for details about viscosity solutions and their properties). As usual, Pucci's classes are defined as
\begin{definition}\label{de.P}
Let $\Omega\subset\mathbb{R}^n$ be a domain and $f\in C(\Omega)$. We write $u\in \underline{S}(\lambda,\Lambda,f)$ in $\Omega$ if $u\in C(\Omega)$ is a viscosity subsolution of
\begin{equation*}
\mathcal{M}^+(D^2u,\lambda,\Lambda)=f\quad\mbox{in}~~\Omega.
\end{equation*}
Similarly, we write $u\in \overline{S}(\lambda,\Lambda,f)$ in $\Omega$ if $u\in C(\Omega)$ is a viscosity supersolution of
\begin{equation*}
\mathcal{M}^-(D^2u,\lambda,\Lambda)=f\quad\mbox{in}~~\Omega.
\end{equation*}
We write $u\in S(\lambda,\Lambda,f)$ in $\Omega$ if $u\in \underline{S}(\lambda,\Lambda,f)$ and $u\in \overline{S}(\lambda,\Lambda,f)$ simultaneously.
\end{definition}

Our main results are the following.
\begin{theorem}\label{t.ful-1}
Let $u\geq 0$ be a viscosity solution of
\begin{equation}\label{e.ful-1}
\left\{\begin{aligned}
&u\in S(\lambda,\Lambda,0)&& ~~\mbox{in}~~ \mathbb{R}^n _+;\\
&u=0&& ~~\mbox{on}~~\partial  \mathbb{R}^n _+.
\end{aligned}\right.
\end{equation}
Then
\begin{equation*}\label{e.ful-1-xn}
  u\equiv u(e_n)x_n\quad\mbox{in}~~\mathbb{R}^n _+.
\end{equation*}
\end{theorem}

The next is a second order Liouville theorem on $ \mathbb{R}^n _+$.
\begin{theorem}\label{t.ful-2}
Let $u\geq 0$ be a viscosity solution of
\begin{equation}\label{e.ful-2}
\left\{\begin{aligned}
&F(D^2u)=1&& ~~\mbox{in}~~ \mathbb{R}^n _+;\\
&u=0&& ~~\mbox{on}~~\partial  \mathbb{R}^n _+.
\end{aligned}\right.
\end{equation}
Then
\begin{equation}\label{e.ful-2-xn2}
  u\equiv ax_n+bx_n^2\quad\mbox{in}~~\mathbb{R}^n _+,
\end{equation}
where $a,b\geq 0$ depend only on $F$.
\end{theorem}

\begin{remark}\label{re1.1}
Braga \cite{MR3735564} proved \Cref{t.ful-1} by comparing the solution with linear functions (i.e., $cx_n$). In this note, we give another proof based on boundary pointwise $C^{1,\alpha}$ regularity although the underlying ideas are essentially the same. The advantage of our proof is that it also applies to \Cref{t.ful-2}, which appears to be new.

It is worth noting that we can only obtain Liouville-type theorems on half-spaces for $S(\lambda,\Lambda,0)$ and for equations of the form $F(D^2u)=1$, without imposing any smoothness or convexity assumptions on the fully nonlinear operators. The reason is that we only have interior $C^{\alpha}$ regularity for $S(\lambda,\Lambda,0)$ and interior $C^{1,\alpha}$ regularity for $F(D^2u)=1$. In contrast, we have boundary $C^{1,\alpha}$ regularity and boundary $C^{2,\alpha}$ regularity, respectively.
\end{remark}

\section{Proofs of main results}\label{sec:1}
Our proofs base on the boundary pointwise regularity heavily. Let us first recall them. The first is the boundary pointwise $C^{1,\alpha}$ regularity (see \cite[Theorem 1.1]{MR3246039}):
\begin{lemma}[\textbf{Boundary pointwise $C^{1,\alpha}$ regularity}]\label{t.reg-S}
Let $u$ be a viscosity solution of
\begin{equation}\label{e.reg-1}
\left\{\begin{aligned}
&u\in S(\lambda,\Lambda,f)&& ~~\mbox{in}~~B_1^+;\\
&u=0&& ~~\mbox{on}~~T_1.
\end{aligned}\right.
\end{equation}
Then $u\in C^{1,\alpha}(0)$, i.e., there exists a constant $a$ such that
\begin{equation}\label{e.reg-res}
    |u(x)-ax_n|\leq C|x|^{1+\alpha}\left(\|u\|_{L^{\infty}(B_1^+)}+\|f\|_{L^{\infty}(B_1^+)}\right), ~~\forall ~x\in B_{1/2}^+
\end{equation}
and
\begin{equation*}
  |a|\leq C\|u\|_{L^{\infty}(B_1^+)},
\end{equation*}
where $0<\alpha<1$ and $C$ are universal (i.e., depending only on $n,\lambda$ and $\Lambda$).
\end{lemma}

The next is the boundary pointwise $C^{2,\alpha}$ regularity (see \cite[Lemma 4.1]{MR3246039}):
\begin{lemma}[\textbf{Boundary pointwise $C^{2,\alpha}$ regularity}]\label{t.reg-S-2}
Let $u$ be a viscosity solution of
\begin{equation}\label{e.reg-2}
\left\{\begin{aligned}
&F(D^2u)=1&& ~~\mbox{in}~~B_1^+;\\
&u=0&& ~~\mbox{on}~~T_1.
\end{aligned}\right.
\end{equation}
Then $u\in C^{2,\alpha}(0)$, i.e., there exist constants $a$ and $b_{in}$($1\leq i\leq n$) such that
\begin{equation}\label{e.reg-res-2}
    |u(x)-ax_n-b_{in}x_ix_n|\leq C|x|^{2+\alpha}\left(\|u\|_{L^{\infty}(B_1^+)}+|F(0)|+1\right), ~~\forall ~x\in B_{1/2}^+
\end{equation}
and
\begin{equation*}
  |a|+|b_{in}|\leq C\left(\|u\|_{L^{\infty}(B_1^+)}+|F(0)|+1\right),
\end{equation*}
where $0<\alpha<1$ and $C$ are universal. Note that the summation over $i$ is understood for the term $b_{in}x_ix_n$ in \cref{e.reg-res-2}.
\end{lemma}

From the boundary pointwise regularity, we can obtain the following Liouville theorems immediately:
\begin{theorem}\label{co2.2}
Let $u$ be a viscosity solution of
\begin{equation}\label{e.reg-2}
\left\{\begin{aligned}
&u\in S(\lambda,\Lambda,0)&& ~~\mbox{in}~~\mathbb{R}_+^n;\\
&u=0&& ~~\mbox{on}~~\partial \mathbb{R}_+^n.
\end{aligned}\right.
\end{equation}
Suppose that for some positive constant $K$,
\begin{equation}\label{e3.2}
|u(x)|\leq K(1+|x|),~\forall ~x\in \mathbb{R}^n_+.
\end{equation}
Then
\begin{equation*}
  u\equiv u(e_n)x_n~~\mbox{ in}~~\mathbb{R}_+^n,
\end{equation*}
\end{theorem}
\proof A scaling version of \Cref{t.reg-S} gives that for any $R>0$,
\begin{equation}\label{e.t.ful-reg-2}
|u(x)-ax_{n}|\leq C\frac{|x|^{1+\alpha}}{R^{1+\alpha}}\cdot \|u\|_{L^{\infty}(B_R^+)},
~~\forall~x\in B^+_{R/2}.
\end{equation}
By the assumption \cref{e3.2},
\begin{equation*}
\|u\|_{L^{\infty}(B_R^+)}\leq K(1+R).
\end{equation*}
Fix $x\in  \mathbb{R}^n _+$ and let $R\rightarrow \infty$ in \cref{e.t.ful-reg-0} (note that $a$ is independent of $R$). Then $u(x)= ax_n$  and hence
\begin{equation*}
  u\equiv ax_n~\mbox{ in }~ \mathbb{R}^n _+.
\end{equation*}
Obviously, $a=u(e_n)$.~\qed

\medskip

\begin{theorem}\label{co2.1}
Let $u$ be a viscosity solution of
\begin{equation}\label{e.reg-1}
\left\{\begin{aligned}
&F(D^2u)=1&& ~~\mbox{in}~~\mathbb{R}_+^n;\\
&u=0&& ~~\mbox{on}~~\partial \mathbb{R}_+^n.
\end{aligned}\right.
\end{equation}
Suppose that for some positive constant $K$,
\begin{equation}\label{e3.1}
|u(x)|\leq K(1+|x|^2),~\forall ~x\in \mathbb{R}^n_+.
\end{equation}
Then
\begin{equation*}
  u\equiv ax_n+b_{in}x_ix_n~~\mbox{ in}~~\mathbb{R}_+^n
\end{equation*}
for some constants $a$ and $b_{in}$.
\end{theorem}
\proof A scaling version of \Cref{t.reg-S-2} gives that for any $R>0$,
\begin{equation}\label{e.t.ful-reg-0}
|u(x)-ax_{n}-b_{in}x_ix_n|\leq C\frac{|x|^{2+\alpha}}{R^{2+\alpha}}\cdot \left(\|u\|_{L^{\infty}(B_R^+)}+R^2|F(0)|+R^2\right),
~~\forall~x\in B^+_{R/2}.
\end{equation}
By the assumption \cref{e3.1},
\begin{equation*}
\|u\|_{L^{\infty}(B_R^+)}\leq K(1+R^2).
\end{equation*}
Fix $x\in  \mathbb{R}^n _+$ and let $R\rightarrow \infty$ in \cref{e.t.ful-reg-0}. Then $u(x)=ax_n+b_{in}x_ix_n$, i.e.,
\begin{equation*}
  u\equiv ax_n+b_{in}x_ix_n~\mbox{ in }~ \mathbb{R}^n _+.
\end{equation*}
~\qed~

\medskip

In the following, we will prove Liouville theorems with the growth conditions \cref{e3.2} and \cref{e3.1} replaced by $u\geq 0$. We will show that the later implies the former. To do so, we need two basic estimates: one is the Carleson type estimate, which was first obtained by Carleson for harmonic functions \cite{MR159013}.
\begin{lemma}[\textbf{Carleson type estimate}]\label{l.cal}
Let $u\geq 0$ be a viscosity solution of
\begin{equation*}\label{e.cal-1}
\left\{\begin{aligned}
&u\in S(\lambda,\Lambda,f)&& ~~\mbox{in}~~B_1^+;\\
&u=0&& ~~\mbox{on}~~T_1.
\end{aligned}\right.
\end{equation*}
Then
\begin{equation}\label{e.cal-2}
  \|u\|_{L^{\infty}(B^+_{1/2})}\leq C\left(u(e_n/4)+\|f\|_{L^{n}(B_1^+)}\right),
\end{equation}
where $C$ is universal.
\end{lemma}

\proof Up to a normalization, we assume
\begin{equation*}
u(e_n/4)+\|f\|_{L^{n}(B_1^+)}\leq 1.
\end{equation*}
Let $v$ be the zero extension of $u$ to the whole $B_1$. Then $v$ satisfies (see \cite[Proposition 2.8]{MR1351007})
\begin{equation*}
 v\in \underline{S}(\lambda,\Lambda,\min(f,0))~~\mbox{ in}~~B_1.
\end{equation*}
By the local maximum principle (see \cite[Theorem 4.8]{MR1351007}), for any $p>0$,
\begin{equation*}\label{cal-lp}
  \|u\|_{L^{\infty}(B_{1/2}^+)}\leq \|v\|_{L^{\infty}(B_{1/2})}
  \leq C_0 \left(\|v\|_{L^{p}(B_{3/4})}+\|f\|_{L^{n}(B_1^+)}\right)= C_0(\|u\|_{L^{p}(B^+_{3/4})}+1),
\end{equation*}
where $C_0$ depends only on $n,\lambda,\Lambda$ and $p$. Thus, we only need to prove that for some universal constants $p>0$ and $C>0$,
\begin{equation*}\label{e2.1}
\|u\|_{L^{p}(B^+_{3/4})}\leq C.
\end{equation*}

Fix $x_0 \in \bar{B}^+_{3/4}$ and denote
\begin{equation*}
L(r,R)=\left\{(x',x_n): x'=x'_0,r\leq x_n\leq R\right\}\cap B^+_{3/4},~\forall ~0<r\leq R.
\end{equation*}
From the Harnack inequality (see \cite[Theorem 4.3]{MR1351007}),
\begin{equation*}
\sup_{L(1/16,3/4)} u\leq C\left(u(e_n/2)+\|f\|_{L^{n}(B_1^+)}\right)\leq C.
\end{equation*}
By considering the equation in
\begin{equation*}
\left\{(x',x_n): |x'-x'_0|<1/32,~~1/64<x_n<1/8\right\}
\end{equation*}
and applying the Harnack inequality again,
\begin{equation*}
  \sup_{L(1/32,1/16)} u\leq C\left(u(x'_0,1/16)+\|f\|_{L^{n}(B_1^+)}\right)\leq C_1,
\end{equation*}
where $C_1$ is universal. Similarly, we have
\begin{equation*}
  \sup_{L(1/2^k,1/2^{k-1})} u\leq C_1^{k-4},~~~~\forall~~k\geq 5.
\end{equation*}

If $x_{0,n}>1/16$ (i.e., the $n$-th component of $x_0$), we have $u(x_0)\leq C$ by the Harnack inequality. Otherwise, there exists $k\geq 5$ such that $1/2^k\leq x_{0,n}\leq 1/2^{k-1}$. Then by taking $q>0$ such that $2^q=C_1$,
\begin{equation}\label{e2.5}
  u(x_0)\leq \sup_{L(1/2^k,1/2^{k-1})} u \leq C_1^{k-5}\leq (2^{-(k-1)})^{-q}\leq x_{0,n}^{-q}.
\end{equation}
Hence,
\begin{equation*}
\|u\|_{L^{\frac{1}{2q}}(B^+_{3/4})}\leq C.
\end{equation*}
\qed \medskip

\begin{remark}\label{re2.2}
The idea that obtain the Carleson type estimate by combining the interior Harnack inequality and zero extension is inspired by the work of De Silva and Savin \cite{MR4093736} (see Step 2 on Page 2423).
\end{remark}

\medskip

The second key estimate is an extension of the Hopf lemma.
\begin{lemma}[\textbf{Hopf type estimate}]\label{le2.2}
Let $u\geq 0$ be a viscosity solution of
\begin{equation*}
\left\{\begin{aligned}
&u\in S(\lambda,\Lambda,f)&& ~~\mbox{in}~~B_1^+;\\
&u=0&& ~~\mbox{on}~~T_1.
\end{aligned}\right.
\end{equation*}
Then
\begin{equation}\label{e3.3}
  u(e_n/4)\leq C\left(u(re_n)/r+\|f\|_{L^{\infty}(B_1^+)}\right),~\forall ~0<r<1/4,
\end{equation}
where $C$ is universal.
\end{lemma}

\medskip

\begin{remark}\label{re2.1}
If $f\equiv 0$, by letting $r\to 0$ in \cref{e3.3}, we have
\begin{equation*}
u_n(0)\geq Cu(e_n/4)>0,
\end{equation*}
which is exactly the classical Hopf lemma. In addition, from \cref{e3.3}, we infer that if $\|f\|_{L^{\infty}(\infty)}\leq c u(e_n/4)$ for some universal small constant $c<1$, the Hopf lemma holds as well (see also \cite{DLL_2025} for a similar result for parabolic equations). Note that this smallness condition is essential. This can be illustrated by the counterexample $u=x_n^2$, which is a nonnegative solution of $\Delta u=1$ in $B_1^+$ and $u=0$ on $T_1$.
\end{remark}

\begin{remark}\label{re2.4}
From the following proof, we know that $f\in L^{\infty}$ can be weakened to $f\in L^p$ for some $p>n$ since the boundary $C^{1,\alpha}$ regularity holds under the later condition as well.
\end{remark}

\medskip

\proof We prove the lemma by contradiction. Suppose the lemma is false. Then there exist sequences of $u_k,f_k,r_k$ satisfying $0<r_k<1/4$ and
\begin{equation*}
\left\{\begin{aligned}
&u_k\in S(\lambda,\Lambda,f_k)&& ~~\mbox{in}~~B_1^+;\\
&u_k=0&& ~~\mbox{on}~~T_1.
\end{aligned}\right.
\end{equation*}
Moreover,
\begin{equation*}
  u_k(e_n/4)\geq k\left(u(r_ke_n)/r_k+\|f_k\|_{L^{\infty}(B_1^+)}\right),~\forall ~k\geq 1.
\end{equation*}
Up to a normalization, we assume
\begin{equation}\label{e3.4}
  u_k(e_n/4)=1, ~\forall ~k\geq 1.
\end{equation}
Hence,
\begin{equation}\label{e3.5}
u(r_ke_n)/r_k+\|f_k\|_{L^{\infty}(B_1^+)}\leq 1/k,~\forall ~k\geq 1.
\end{equation}

By the Carleson type estimate (see \Cref{l.cal}),
\begin{equation*}
  \|u_k\|_{L^{\infty}(B_{1/2}^+)}\leq C\left(u_k(e_n/4)+\|f_k\|_{L^{n}(B_1^+)}\right)\leq C.
\end{equation*}
From the global $C^{\alpha}$ regularity,
\begin{equation*}
  \|u_k\|_{C^{\alpha}(\bar B_{3/8}^+)}\leq C\left(\|u_k\|_{L^{\infty}(B_{1/2}^+)}+\|f_k\|_{L^{n}(B_1^+)}\right)\leq C.
\end{equation*}
Hence, up to a subsequence, there exists $\bar{u}\in C(\bar{B}_{3/8}^+)$ such that $u_k\to \bar{u}$ uniformly on $\bar{B}_{3/8}^+$. Then $\bar{u}\geq 0$ is a viscosity solution of
\begin{equation*}
\left\{\begin{aligned}
&\bar u\in S(\lambda,\Lambda,0)&& ~~\mbox{in}~~B_1^+;\\
&\bar u=0&& ~~\mbox{on}~~T_1.
\end{aligned}\right.
\end{equation*}
Furthermore, by \cref{e3.4},
\begin{equation}\label{e3.6}
\bar{u}(e_n/4)=1.
\end{equation}

Up to a subsequence, we have $r_k\to \bar{r}$ for some $\bar{r}\geq 0$. If $\bar{r}>0$, by letting $k\to \infty$ in \cref{e3.5}, we have
\begin{equation*}
  \bar{u}(\bar{r}e_n)=0,
\end{equation*}
which is a contradiction to the strong maximum principle.

Next, we consider the case $\bar{r}=0$. By the boundary $C^{1,\alpha}$ regularity for $u_k$, there exist constants $a_k\geq 0$ such that
\begin{equation}\label{e3.7}
    |u_k(x)-a_kx_n|\leq C|x|^{1+\alpha}\left(\|u_k\|_{L^{\infty}(B_{1/2}^+)}+\|f_k\|_{L^{\infty}(B_1^+)}\right)
    \leq C|x|^{1+\alpha}, ~~\forall ~x\in B_{3/8}^+.
\end{equation}
Choose $x=r_ke_n$ in above inequality. Then we have
\begin{equation*}
|u_k(r_ke_n)/r_k-a_k|\leq Cr_k^{\alpha}.
\end{equation*}
Let $k\to \infty$ and note \cref{e3.5}. Then
\begin{equation*}
  a_k\to 0 \quad\mbox{as}~~k\to \infty.
\end{equation*}
Finally, take $k\to \infty$ in \cref{e3.7} and we have
\begin{equation*}
|\bar{u}(x)|\leq C|x|^{1+\alpha}, ~~\forall ~x\in B_{3/8}^+,
\end{equation*}
which implies
\begin{equation*}
\bar{u}_n(0)=0.
\end{equation*}
This is impossible by the Hopf lemma.~\qed~

\medskip

Now, we can prove the Liouville theorems on half-spaces. First, we give the~\\
\noindent\textbf{Proof of \Cref{t.ful-1}.} For any $R>0$, from the Hopf type estimate (\Cref{le2.2}),
\begin{equation}\label{e2.3}
  u(Re_n/4)\leq C u(e_n/2) R.
\end{equation}
In addition, the Carleson type estimate (\Cref{l.cal}) implies
\begin{equation}\label{e2.2}
  \|u\|_{L^{\infty}(B_R^+)}\leq Cu(Re_n/4).
\end{equation}
Hence,
\begin{equation*}
 \|u\|_{L^{\infty}(B_R^+)}\leq C u(e_n/2) R.
\end{equation*}
That is, $u$ grows linearly at most. By \Cref{co2.2},
\begin{equation*}
  u\equiv u(e_n)x_n~~~~\mbox{  in }~ \mathbb{R}^n _+.
\end{equation*}
~\qed~

\medskip

\begin{remark}\label{r-1.2}
The above proof shows the idea clearly. The boundary pointwise regularity gives the estimate of the error between the solution and a linear function (see \Cref{t.reg-S}). If $\|u\|_{L^{\infty}(B_R^+)}$ grows linearly as $R\rightarrow \infty$, we arrive at the Liouville theorem (see \Cref{co2.2}). The Hopf type lemma guarantees that $u(Re_n/2)$ grows linearly at most (see \cref{e2.3}). Finally, the Carleson type estimate provides a bridge between $\|u\|_{L^{\infty}(B_R^+)}$ and $u(Re_n/2)$ (see \cref{e2.2}).
\end{remark}

\medskip

With the same idea, we can prove the second order Liouville theorem:~\\
\noindent\textbf{Proof of \Cref{t.ful-2}.} For any $R>0$, from the Hopf type estimate (\Cref{le2.2}),
\begin{equation}\label{e2.3-2}
  u(Re_n/4)\leq C \left(u(e_n/2) R+R^2|F(0)|+R^2\right).
\end{equation}
Moreover, the Carleson type estimate (\Cref{l.cal}) shows
\begin{equation}\label{e2.2-2}
  \|u\|_{L^{\infty}(B_R^+)}\leq C\left(u(Re_n/4)+R^2|F(0)|+R^2\right).
\end{equation}
Thus,
\begin{equation*}
 \|u\|_{L^{\infty}(B_R^+)}\leq C\left( u(e_n/2) R+R^2\right).
\end{equation*}
That is, $u$ grows quadratically at most. By \Cref{co2.1} and noting $u\geq 0$, we have
\begin{equation*}
  u\equiv ax_n+bx^2_n \quad\mbox{in }~ \mathbb{R}^n _+
\end{equation*}
for some constants $a,b\geq 0$. ~\qed~\\

\printbibliography

\end{document}